\documentclass[10 pt]{article}
\usepackage[centertags]{amsmath}
\usepackage{amssymb}
\usepackage{amsthm}
\usepackage{newlfont}
\usepackage{amsfonts}

\usepackage{graphicx}
\usepackage{color}
\usepackage[colorlinks]{hyperref}
\usepackage{mathptmx}

\newtheorem{thm}{Theorem}

\newcommand{\bks}{\backslash}
\newcommand{\tbf}{\textbf}

\begin{document}

\begin{center}
  New Proofs of K\"{o}nig's Bipartite Graph Characterization Theorem
\end{center}

\begin{center}
  Salman Ghazal\footnote{Department of Mathematics, Faculty of Sciences I, Lebanese University, Hadath, Beirut, Lebanon.\\
                       E-mail: salmanghazal@hotmail.com\\
                       }
\end{center}

\begin{abstract}
We introduce four new elementary short proofs of the famous K\"{o}nig's theorem which characterizes bipartite graphs by absence of odd cycles.

\end{abstract}

\section{Introduction}

In this short paper, graphs are finite and may contain loops or multiple edges. The vertex set of a graph $G$ is denoted by $V(G)$ while its edge set is denoted by $E(G)$. The induced subgraph of $G$ by $A\subseteq V(G)$ is denoted by $G[A]$. A subgraph $H$ of $G$ is called a spanning subgraph of $G$ if $V(H)=V(G)$. A set $X$ of pairwise nonadjacent vertices of $G$ is said to be a \emph{stable} set of $G$, that is $G[X]$ has no edge. With abuse of notation, $xy$ is used to denoted an edge whose endpoints are the vertices $x$ and $y$. The length of a path or cycle is the number of its edges. A cycle of odd (resp. even) length is called an odd (resp. even) cycle. A path between two vertices $a$ and $b$ is called an $a,b$-path. We do not distinguish between a connected component and the subgraph it induces.\\

A graph $G$ is \emph{bipartite} if its vertex set is the union of two disjoint (possibly empty) stable sets $X$ and $Y$. In this case, $\{X,Y\}$ is said to be a \emph{bipartition} of $G$.\\

Let $G$ be a graph. It is clear that $G$ is bipartite if and only if all its connected components are so. Moreover, If $G'$ is obtained from $G$ by keeping only one copy of each set multiple edges of $G$, then $G$ is bipartite if and only if $G'$ is so. In addition, if $\{X, Y\}$ is a bipartition of $G$, $a\in X$ and $b\in Y$, then $\{X,Y\}$ is again a bipartition of $G+ab$, because $X$ and $Y$ are still stable in $G+ab$.\\

In fact, suppose that $A_1,\cdots , A_k$ are the connected components of a bipartite graph $G$ with bipartition $\{X, Y\}$. For $i= 1, \cdots, k$, let $X_i=X\cap A_i$ and $Y_i=Y\cap A_i$. Then $\{X_i,Y_i\}$ is a bipartition of the connected component $A_i$. Moreover, $\forall 1\leq i\leq k$, the sets $X'=(X-X_i)\cup Y_i$ and $Y'=(Y-Y_i)\cup X_i$ form a bipartition of $G$.\\

Suppose that $P=x_1x_2...x_n$ is a path in a bipartite graph $G$ with a specified bipartition $\{X,Y\}$. Note that if $a$ and $b$ are adjacent vertices of $G$, then they must be in distinct partite sets. So, if $x_1\in X$, then so is every vertex of $P$ with odd index, while every vertex of $P$ with even index is in $Y$. Hence, $n$ is odd if and only if $x_n\in X$. Therefore, if $C=x_1x_2\cdots x_nx_1$ is a cycle of $G$, then it must be even, since otherwise the adjacent vertices $x_n$ and $x_1$ must be in the same partite set, which contradicts its stability.\\

In fact, the above obvious necessary condition of bipartite graphs is also sufficient. This is proved in 1936 by K\"{o}nig \cite{konig}. Proofs of the sufficient condition used distances, walks or spanning trees.

\begin{thm} \tbf{(K\"{o}nig \cite{konig})}
A graph is bipartite if and only if it has no odd cycle.
\end{thm}

\section{Four New Elementary Proofs}

We introduce three new elementary proofs of the sufficient condition of K\"{o}nig's theorem that use neither distances nor walks nor spanning trees. We may assume that $G$ has no multiple edges.\\

\noindent\textbf{First Proof:}
\begin{proof}
Let $G$ be a graph that has no odd cycle. We may assume that $G$ is connected. Since $G$ has no loop, any vertex of $G$ can be viewed as a bipartite, connected and induced subgraph of $G$. Let $H$ be a maximal bipartite, connected and induced subgraph of $G$. We prove that $G=H$ and consequently we get that $G$ is a bipartite graph. Suppose to the contrary that $G\neq H$. Then $V(H)\neq V(G)$. Since $G$ is connected, $ \exists z\in V(G)\bks V(H)$ and $\exists t\in V(H)$ such that $zt\in E(G)$. Let $\{X_1 ,X_2\}$ be a bipartition of $H$. For $i=1,2$, if $\forall x\in X_i$, $zx\notin E(G)$, then $X_{i}\cup \{z\}$ would be a stable set and thus $G[V(H)\cup \{z\}]$ would be a bipartite, connected, induced subgraph of $G$ and that contains $H$ strictly, which contradicts the maximality of $H$. Hence, for $i=1,2$, $\exists x_i\in X_i$ such that $zx_i\in E(G)$. However, since $H$ is connected, it contains an $x_1,x_2$-path $P$. Since $H$ is bipartite and $x_1$ and $x_2$ are in distinct partite sets, then the length of $P$ is odd. Therefore, adding to $P$ the edges $zx_1$ and $zx_2$ forms an odd cycle, which is a contradiction.

\end{proof}

\noindent\textbf{Second Proof:}

\begin{proof}
Let $G$ be a graph that has no odd cycle. The spanning subgraph of $G$ with no edges is bipartite. Let $H$ be a maximal bipartite spanning subgraph of $G$. We prove that $G=H$ and consequently we get that $G$ is a bipartite graph. Suppose to the contrary that $G\neq H$. Then $E(H)\neq E(G)$ and hence $\exists e=ab\in E(G)-E(H)$. Let $\{X ,Y\}$ be a bipartition of $H$. By maximality of $H$, the graph $H'=H+e$ is not bipartite and thus $a$ and $b$ lie in the same partite set of $H$, say $X_1$, since otherwise, $\{X ,Y\}$ would be a bipartition of $H'$ also. If there is an $ab$-path $P$ in $H$, then its length is even and adding to it the edge $e$ would create an odd cycle in $G$, a contradiction. Therefore, $a$ and $b$ are in distinct components of $H$. Let $A$ be the connected component of $H$ containing $a$. Then $X'=(X-(X\cap A))\cup (Y\cap A)$ and $Y'=(Y-(Y\cap A))\cup (X\cap A)$ is a bipartition of $H$ and $H'$. This contradicts the fact that $H'$ is not bipartite.
\end{proof}

\noindent\textbf{Third Proof:}
\begin{proof}
 Suppose that a counterexample exist and let $G$ be a minimal one. Let $e=ab\in E(G)$. Then $G-e\subsetneq G$ and $G-e$ has no odd cycle. Hence, $G-e$ is bipartite, by minimality of $G$. Let $\{X, Y\}$ be a bipartition of $G-e$. If $a$ and $b$ are in distinct partite sets, then $\{X, Y\}$ is a bipartition of $G$ as well, a contradiction. So, $a$ and $b$ belong to the same partite set, say $X$.
 Suppose that there is an $a,b$-path $P\subseteq G$ distinct from $ab$. Then $P\subseteq G-e$. Since $a$ and $b$ are in the same partite set of $G-e$, then the length of $P$ even. Adding to $P$ the edge $ab$ creates an odd cycle in $G$, a contradiction. So $ab$ is the unique $a,b$-path. Thus $G-e$ is not connected. Let $A$ be the connected component of $G-e$ containing $a$. Then $b\notin A$. Now $X'=(X-(X\cap A))\cup (Y\cap A)$ and $Y'=(Y-(Y\cap A))\cup (X\cap A)$ form a bipartition of $G-e$ and thus $G$, since  $b\in X'$ and $a\in Y'$. %Hence, $\{X', Y'\}$ is a bipartition of $G$.
  A contradiction.

 \end{proof}
 
 \noindent\textbf{Fourth Proof:}
\begin{proof}

First we prove by induction on the number of vertices that if a graph has no cycle, then it is bipartite. Let $F$ be such a graph. Then $F$ has a vertex $x$ that has at most one neighbor $y$. Since $F-x$ has no cycles as well, then by the induction hypothesis, it is bipartite, with bipartition say $\{A,B\}$. We may assume that $y\notin A$ (if $y$ exist). Then $\{A\cup\{x\}, B\}$ is a bipartition of $F$.\\

Let $G$ be a graph that has no odd cycle. If $G$ has no (even) cycle, then it is bipartite. Otherwise let $e\in E(C)$, for some even cycle $C$ of $G$. By induction on the number of cycles of $G$, the graph $G-e$ is bipartite since it has fewer cycles than $G$. But the path $P=C-e\subseteq G-e$ has an odd length, hence its endpoints are in distinct partite sets. Thus $G=(G-e)+e$ is bipartite.

\end{proof}

\end{document}